\documentclass{article}
\usepackage{amsmath,amssymb,amscd}

\newtheorem{theorem}{Theorem}[section]

\newtheorem{corollary}[theorem]{Corollary}

\newtheorem{question}[theorem]{Question}

\newcommand{\pri}{\ensuremath{\smallsetminus}}

\def\nl{[}
\def\np{]}

\makeatletter
\def\@cite#1#2{\nl{#1\if@tempswa , #2\fi}\np}
\makeatother

\begin{document}

\begin{center}
 {\LARGE\bf Closed subgroups of the polynomial automorphism group
  containing the affine subgroup.}\\
\vspace{.4cm}
{\large Eric Edo}\\
\end{center}
\vspace{.4cm}

\begin{abstract}
We prove that, in characteristic zero, 
closed subgroups of the polynomial automorphisms group containing the affine group contain the whole tame group.
\end{abstract}

\section{Introduction}

\ \ \ Throughout, $\mathbb{K}$ denotes an algebraically closed field of characteristic ${\rm char}(\mathbb{K})$
and $n\in\mathbb{N}_+$ denotes a positive integer.\\

We denote by ${\cal E}={\rm End}_n(\mathbb{K})$ the set $(\mathbb{K}[x_1,\ldots,x_n])^n$ endowed
with the following law
$\sigma\tau=(\,f_1(g_1,\ldots,g_n)\,,\ldots,\,f_n(g_1,\ldots,g_n)\,)$
for all elements $\sigma=(f_1,\ldots,f_n)$ and $\tau=(g_1,\ldots,g_n)$ in~${\cal E}$.
The monoid ${\cal E}$ is anti-isomorphic to ${\rm End}_{\mathbb{K}}\mathbb{K}[x_1,\ldots,x_n]$ 
(the monoid of polynomial endomorphisms).\\

We denote by ${\cal G}={\rm GA}_n(\mathbb{K})$ the group of invertible elements of ${\cal E}$.
Note that ${\cal G}$ is anti-isomorphic to the group
${\rm Aut}_{\mathbb{K}}\mathbb{K}[x_1,\ldots,x_n]$ (the group of polynomial automorphisms)
and isomorphic to the group of polynomial automorphisms of the affine $n$-space $\mathbb{A}^n_{\mathbb{K}}$.\\

We define the degree of $\sigma=(f_1,\ldots,f_n)\in{\cal E}$ by $\deg(\sigma)=\max_{1\le i\le n}\{\deg(f_i)\}$.
We denote by ${\cal A}={\rm Aff}_n(\mathbb{K})=\{\sigma\in{\cal G}\,;\,\deg(\sigma)=1\}$ the \textit{affine subgroup}, by
${\cal B}={\rm BA}_n(\mathbb{K})= 
\{(\alpha_1x_1+p_1,\ldots,\alpha_nx_n+p_n);
\alpha_i\in\mathbb{K}^*,p_i\in k[x_{i+1},\ldots,x_n],\forall i\}$
the \textit{triangular subgroup}
and by ${\cal T}={\rm TA}_n(\mathbb{K})$ the \textit{tame subgroup} i.~e. the subgroup 
generated by ${\cal A}$ and ${\cal B}$
 (see \cite{vdE} for more informations about polynomial automorphisms).\\

Given a positive integer $d\in\mathbb{N}_+$, the set ${\cal E}_{\le d}$ of polynomial endomorphisms of degree $\le d$ is
a $\mathbb{K}$-vector space of dimension $N={n+d\choose d}$ (the degree of a polynomial endomorphism is the maximum of the degree of its components)
and we can consider ${\cal E}_{\le d}$ as an algebraic variety by transfer of the structure of the affine $N$-space $\mathbb{A}^N_{\mathbb{K}}$.
We consider the Zariski topology on ${\cal E}_{\le d}$ associated with this structure.\\

Following Shafarevich (see~\cite{S}), since for all positive integers $d\in\mathbb{N}_+$, ${\cal E}_{\le d}$ is closed in ${\cal E}_{\le d+1}$,
we can endow ${\cal E}=\bigcup_{d\ge 1}{\cal E}_{\le d}$ with the inductive limit topology. A subset $A\subset{\cal E}$ is closed in ${\cal E}$
for this topology if and only if $A\cap{\cal E}_{\le d}$ is closed in ${\cal E}_{\le d}$ for all $d\in\mathbb{N}_+$. We consider the restriction
of this inductive limit topology on ${\cal G}$ and we set ${\cal G}_{\le d}={\cal G}\cap{\cal E}_{\le d}$.\\

In this paper, we are interested with the question:

\begin{question}
Which are the closed subgroups of ${\cal G}$ containing ${\cal A}$?
\end{question}

We prove the following result in the characteristic zero case:

\begin{theorem}[$n\ge 2$, ${\rm char}(\mathbb{K})=0$]\label{thm:mainz}
If ${\cal H}$ is a closed subgroup of ${\cal G}$ strictly containing ${\cal A}$ then ${\cal T}\subset{\cal H}$.
\end{theorem}

In~\cite{B1} (see Theorem~2.10), Bodnarchuk claimed to prove this result with even the stronger conclusion
${\cal H}={\cal G}$. But his proof is based on a erroneous result of Shafarevich (cited as Theorem~2.7 in~\cite{B1}).
Recently, Furter and Kraft (see~\cite{FK}) proved that, in the case $n=3$ and $\mathbb{K}=\mathbb{C}$, the subgroup ${\cal T}'$
is closed in ${\cal G}'$ where ${\cal G}'$ is the subgroup of ${\cal G}$ fixing $x_3$ and ${\cal T}'={\cal G}'\cap{\cal T}$.
This result is based on the Shestakov-Umirbaev theory (see~\cite{SU}) which implies that ${\cal T}'\ne{\cal G}'$.
With notations of Theorem~2.7 in~\cite{B1}, taking $G={\cal G}'$, $H={\cal T}'$ and $f:{\cal T}'\to{\cal G}'$ the canonical inclusion,
$({\rm d}f)_e$ is the canonical isomorphism between the tangent spaces $T_{e,H}$ and $T_{e,G}$, but $f$ is not
an isomorphism because of the Shestakov-Umirbaev theorem.\\

In the case $n=2$, the Jung-van der Kulk theorem (see~\cite{J,vdK}) says that ${\cal T}={\cal G}$ and Theorem~\ref{thm:mainz} implies that the affine group is a maximal closed subgroup of the plane polynomial
automorphism group. This
gives a positive answer to a question of Furter (see Question~1 and Question 1.11 in~\cite{F}).\\

 In the case $n=3$, we know that ${\cal T}\ne{\cal G}$.
Recently, Poloni and the author (see~\cite{EP}) proved that some explicit families of automorphisms
in ${\cal G}\pri{\cal T}$ are in the closure of ${\cal T}$ in ${\cal G}$ (see also~\cite{K2}). But we don't know if ${\cal T}$ is dense in ${\cal G}$ or not.
In particular, we don't know if the Nagata automorphism (see~\cite{N}) is in the closure of ${\cal T}$ in ${\cal G}$.\\

The situation in positive characteristic is more complex and will be studied in an other paper.

\section{Closed subgroup}

\ \ \ In this section, we use notations introduced in the last section.
We do not assume that $\mathbb{K}$ has characteristic zero.
We consider the following two sets of variables ${\bf x}=\{x_1,\ldots,x_n\}$ and $\hat{\bf x}=\{x_2,\ldots,x_n\}$.
We denote by $I$ (resp.~$\hat{I}$) the ideal of $\mathbb{K}[{\bf x}]=\mathbb{K}[x_1,\ldots,x_n]$ 
(resp.~$\mathbb{K}[\hat{\bf x}]=\mathbb{K}[x_2,\ldots,x_n]$) generated by ${\bf x}$ (resp.~$\hat{\bf x}$).
We say that an automorphism $\phi=(f_1,\ldots,f_n)\in{\cal G}$ has {\it a affine part equal to} ${\rm id}$
if $f_i\equiv x_i$ modulo $I^2$ for all $i\in\{1,\ldots,n\}$.
We prove the following result which is the main theorem of this paper.

\begin{theorem}[$n\ge 2$]\label{thm:closed}
Let $\phi=(f_1,\ldots,f_n)\in{\cal G}\pri{\cal A}$ be an automorphism of degree $d\ge 2$.
We assume that $\phi$ has a affine part equal to ${\rm id}$ and $f_1\ne x_1$.\\
a) The polynomial $g_0(\hat{\bf x}):=f_1(0,\hat{\bf x})$ is not zero and has a $\hat{I}$-adic valuation
$w\ge 2$.\\
b) Let ${\cal H}:=\langle{\cal A},\phi\rangle$ be the subgroup of ${\cal G}$
generated by ${\cal A}\cup\{\phi\}$. The closure of ${\cal H}\cap{\cal G}_{\le d}$ in ${\cal G}_{\le d}$
contains the automorphism $(x_1+h(\hat{\bf x}),\hat{\bf x})$
where $h(\hat{\bf x})$ is the homogeneous part of $g_0(\hat{\bf x})$ of degree $w$.
\end{theorem}

{\bf Proof\ }
We write $f_1=\sum_{k=0}^mg_kx_1^k$ where $m=\deg_{x_1}f_1\ge 1$ and $g_i\in\mathbb{K}[\hat{\bf x}]$.\\
a) By contradiction, we assume that $g_0=0$. Then $x_1$ divides $f_1$. Since $f_1$ is a coordinate, $f_1$ is an irreducible polynomial and using that $f_1\equiv x_1$ modulo $I^2$, we deduce $f_1=x_1$. Impossible. 
Since the affine part of $\phi$ is ${\rm id}$,
it's clear that $2\le w<+\infty$.\\
b) We consider the torus action $\alpha_t:=(t^wx_1,t\hat{\bf x})\in{\rm Aff}_n(\mathbb{K}(t))$ 
and we compute
$$\phi_t:=\alpha_t^{-1}\phi\alpha_t
=(t^{-w}f_1(t^wx_1,t\hat{\bf x}),t^{-1}f_2(t^wx_1,t\hat{\bf x}),\ldots,t^{-1}f_n(t^wx_1,t\hat{\bf x})).$$
Using that the affine part of $\phi$ is ${\rm id}$, we have $g_1(0)=1$ and
$$t^{-w}f_1(t^wx_1,t\hat{\bf x})= t^{-w}\sum_{k=0}^mg_k(t\hat{\bf x})t^{kw}x_1^k
\equiv t^{-w}g_0(t\hat{\bf x})+x_1\equiv x_1+h(\hat{\bf x})\ \ {\rm mod}\ \  t\mathbb{K}[{\bf x},t]$$
and $t^{-1}f_i(t^wx_1,t\hat{\bf x})\equiv x_i\ \ {\rm mod}\ \ t\mathbb{K}[{\bf x},t]$, for all $i\in\{2,\ldots,n\}$.
Using that ${\rm Jac}(\phi_t)={\rm Jac}(\phi)\in\mathbb{K}^*$,
we deduce from the overring principle (see Lemma~1.1.8 p. 5 in~\cite{vdE}) that $\phi_t\in{\rm GA}_n(\mathbb{K}[t])$.
For all $t_0\in\mathbb{K}^*$, the automorphism $\phi_{t\to t_0}$ is in ${\cal H}\cap{\cal G}_{\le d}$ and we
deduce that $(x_1+h(\hat{\bf x}),\hat{\bf x})=(x_1+h(\hat{\bf x}),x_2,\ldots,x_n)=\phi_{t\to 0}$ is in  the closure of ${\cal H}\cap{\cal G}_{\le d}$ 
in ${\cal G}_{\le d}$.\\

\textbf{Remark:} This proof is based on a classical technique which consists to conjugate an automorphism by an action of the torus.
See the proof of Theorem~4.6 in \cite{EP} or the proof of Lemma~4.1 in \cite{BC}
for other situations where this technique is used.

\begin{corollary}\label{coro:tri}
Let $\phi\in{\cal G}\pri{\cal A}$ be an automorphism of degree $d\ge 2$.
Let ${\cal H}:=\langle{\cal A},\phi\rangle$ be the subgroup of ${\cal G}$
generated by ${\cal A}\cup\{\phi\}$. Then the closure of ${\cal H}\cap{\cal G}_{\le d}$ in ${\cal G}_{\le d}$
contains an element of ${\cal B}\pri{\cal A}$.
\end{corollary}

{\bf Proof\ }
Changing $\phi$ to $\alpha^{-1}\phi$ where $\alpha\in{\cal A}$ is the affine part or $\phi$,
we can assume that the affine part of $\phi$ is ${\rm id}$.
Since $\phi=(f_1,\ldots,f_n)$ is not affine, there exists $i\in\{1,\ldots,n\}$ such that $f_i\ne x_i$.
Changing $\phi$ to $\sigma\phi\sigma$ where $\sigma\in\mathfrak{S}_n\subset{\cal A}$ is the transposition $x_1\leftrightarrow x_i$,
we can assume that $f_1\ne x_1$ and apply Theorem~\ref{thm:closed}.\\

\section{Characteristic zero}

\ \ \ In this section, we assume that $\mathbb{K}$ has characteristic zero.
We recall two complementary results.
Theorem~\ref{thm:F}, in the case $n=2$ (we recall that when $n=2$, we have ${\cal T}={\cal G}$), is due to Furter (see Theorem~D in~\cite{F}).
Theorem~\ref{thm:B}, in the case $n\ge 3$, is due to Bodnarchuk (see ~\cite{B2}).

\begin{theorem}[$n=2$, ${\rm char}(\mathbb{K})=0$]\label{thm:F}
For all $\beta\in{\cal B}\pri{\cal A}$ we have $\overline{\langle{\cal A},\beta\rangle}={\cal G}$.
\end{theorem}

\begin{theorem}[$n\ge 3$, ${\rm char}(\mathbb{K})=0$] \label{thm:B}
For all $\beta\in{\cal B}\pri{\cal A}$ we have $\langle{\cal A},\beta\rangle={\cal T}$.
\end{theorem}

Using this two results, Corollary~\ref{coro:tri} implies Theorem~\ref{thm:mainz}.\\

{\noindent\textbf{Theorem~\ref{thm:mainz}} ($n\ge 2$, ${\rm char}(\mathbb{K})=0$)\textbf{.}
{\it If ${\cal H}$ is a closed subgroup of ${\cal G}$ strictly containing ${\cal A}$ then ${\cal T}\subset{\cal H}$.}}

{\bf Proof\ }
Since ${\cal H}$ strictly containing ${\cal A}$, there exists $\phi\in{\cal H}\pri{\cal A}$ of degree $d\ge 2$.
By Corollary~\ref{coro:tri}, the closure of $\langle{\cal A},\phi\rangle$ in ${\cal G}_{\le d}$ contains an element
$\beta\in{\cal B}\pri{\cal A}$. Using that ${\cal H}$ is closed in ${\cal G}$ we deduce $\beta\in\overline{\langle{\cal A},\phi\rangle}\subset\overline{{\cal H}}={\cal H}$.\\
a) If $n=2$, by Theorem~\ref{thm:F}, we have ${\cal T}={\cal G}=\overline{\langle{\cal A},\beta\rangle}\subset\overline{{\cal H}}={\cal H}$.\\
b)  If $n\ge 3$, by Theorem~\ref{thm:B}, we have ${\cal T}=\langle{\cal A},\beta\rangle\subset{\cal H}$.


\begin{thebibliography}{ASM}

\bibitem{BC} J. Blanc, A. Calabri, {\it On degenerations of plane Cremona transformations},
http://arxiv.org/pdf/1401.6676.pdf

\bibitem{B1} Y. Bodnarchuk, \textit{Some extreme properties of the affine group as an automorphisms group of the affine space},
Contrib. General Algebra, \textbf{13}, (2001) p. 9--22.

\bibitem{B2} Y. Bodnarchuk, {\it Generating properties of triangular and bitriangular birational automorphisms of affine space},
 Dopov. NAN Ukrainy, (2002), no. 11, 7?-22.







\bibitem{EP} E. Edo, P.-M. Poloni, \textit{On the closure of the tame automorphism group of affine three-space},
Inter. Math. Research Notices, (2015), to appear, http://arxiv.org/pdf/1403.2843.pdf

\bibitem{vdE} A. van den Essen, {\it Polynomial Automorphisms
and the Jacobian Conjecture,} Birkhauser Verlag, Basel-Boston-Berlin (2000).



\bibitem{F} J.-P. Furter, {\it Polynomial Composition Rigidity and Plane Polynomial Automorphisms},
http://perso.univ-lr.fr/jpfurter/

\bibitem{FK} J.-P. Furter, H. Kraft, {\it On the geometry of the automorphism group of affine $n$-space},
in preparation, (2015).


\bibitem{J} H.W.E. Jung, \textit{\"{U}ber ganze birationale Transformationen der Ebene},
J. Reine Angew. Math. { 184} (1942), 161-174.

\bibitem{vdK} W. van der Kulk, \textit{On polynomial rings in two variables}, Nieuw. Arch. Wisk. (3) {\bf 1}
(1953), 33-41.


\bibitem{K2} S.~Kuroda, \textit{Degeneration of tame automorphisms of a polynomial ring}, http://arxiv.org/pdf/1411.2085.pdf

\bibitem{N} M.~Nagata, \textit{On Automorphism Group of $k[x,y]$}, Lectures in Mathematics, Department of Mathematics,
Kyoto University, Vol.~5, Kinokuniya Book-Store Co.\ Ltd., Tokyo, 1972.

\bibitem{S} I.R. Shafarevich, {\it On some infinite-dimensional groups II},
Math. USSR Izv, {\bf 18} (1982), 214-226.

\bibitem{SU} I.~Shestakov, U.~Umirbaev, {\it The tame and the wild automorphisms of polynomial rings in three variables},
J.\ Amer.\ Math.\ Soc.\ {\bf 17} (2004), 197--227.


\end{thebibliography}
\end{document}